\documentclass[letterpaper, 12pt]{article}

\usepackage[all,cmtip]{xy}
\usepackage{amsmath}
\usepackage{amssymb}
\usepackage{amsthm}
\usepackage{hyperref}
\usepackage{enumerate}
\usepackage[margin=1in]{geometry}

\newtheorem*{example*}{Example}
\newtheorem*{remark*}{Remark}
\newtheorem*{theorem*}{Theorem}
\newtheorem*{definition*}{Definition}
\newtheorem{theorem}{Theorem}[section]
\newtheorem{rem}[theorem]{Remark}
\newenvironment{remark}{\begin{rem}\rm}{\end{rem}}

\newtheorem{proposition}[theorem]{Proposition}

\newtheorem{lemma}[theorem]{Lemma}

\newtheorem{eg}[theorem]{Example}

\newtheorem{definition}[theorem]{Definition}
\newcommand{\R}{\mathbb{R}}
\renewcommand{\S}{\mathcal{S}}
\newcommand{\C}{\mathbb{C}}
\newcommand{\Z}{\mathbb{Z}}
\newcommand{\W}{\mathcal{W}}
\newcommand{\Cl}{\mathbb{C}\mathrm{l}}
\newcommand{\dirac}{\mathsf{D}_b\!\!\!\!\!\!/\,\,\,}
\newcommand{\dob}{\mathsf{D}\!\!\!\!/\,}

\newcommand{\dbbar}{\overline{\partial}_b}

\DeclareMathOperator{\Lie}{Lie}
\DeclareMathOperator{\End}{End}

\DeclareMathOperator{\Ch}{Ch}

\DeclareMathOperator{\Id}{Id}
\DeclareMathOperator{\im}{im}
\DeclareMathOperator{\Td}{Td}

\DeclareMathOperator{\ind}{index}

\begin{document}
\title{On the geometric quantization of contact manifolds}

\author{Sean Fitzpatrick\\University of California, Berkeley\\sean@math.berkeley.edu}
\maketitle

\begin{abstract}
Suppose that $(M,E)$ is a compact contact manifold, and that a compact Lie group $G$ acts on $M$ transverse to the contact distribution $E$.  In an earlier paper, we defined a $G$-transversally elliptic Dirac operator $\dirac$, constructed using a Hermitian metric $h$ and connection $\nabla$ on the symplectic vector bundle $E\rightarrow M$, whose equivariant index is well-defined as a generalized function on $G$, and gave a formula for its index.  By analogy with the geometric quantization of symplectic manifolds, the $\mathbb{Z}_2$-graded Hilbert space $Q(M)=\ker \dirac \oplus \ker \dirac^{\!\!*}$ can be interpreted as the ``quantization'' of the contact manifold $(M,E)$; the character of the corresponding virtual $G$-representation is then given by the equivariant index of $\dirac$.  By defining contact analogues of the algebra of observables, pre-quantum line bundle and polarization, we further extend the analogy by giving a contact version of the Kostant-Souriau approach to quantization, and discussing the extent to which this approach is reproduced by the index-theoretic method.
\end{abstract}

\section{Introduction}
The problem of geometric quantization is well-known in symplectic geometry, and dates back to the work of Souriau \cite{Sou} and Kostant \cite{Kost}.  Symplectic geometry is the natural setting for classical Hamiltonian dynamics, 
but contact structures appear in classical physics as well: the role of contact geometry in Lagrangian mechanics is explained in \cite{Sev}, and the geometry of classical thermodynamics has a natural contact structure (see  for example \cite{Burke} or \cite{Raj}, which discusses the quantization of thermodynamics via the deformation quantization of contact structures).  Parts of our construction for general contact manifolds reduce to the definitions used in \cite{Raj} when expressed in terms of a local Darboux chart.  Other examples in the literature related to the quantization of contact manifolds include \cite{BMG,GS2,LTW}; in each case the methods used are related to deformation quantization.  A brief sketch of an approach to geometric contact quantization was given by Vaisman in \cite{Vais1}; the first quantization we present for contact manifolds expands upon the suggestion in \cite{Vais1}.  We should also note that a geometric quantization for Jacobi manifolds has been given in \cite{LMP} which specializes to contact manifolds.  However, this approach is based on Vaisman's method of contravariant derivatives in Poisson geometry \cite{Vais2}, while we make use of covariant derivatives, as is the norm in symplectic geometry.

In this article, we will instead describe two ways to define a ``geometric quantization'' of contact manifolds analogous to familiar methods in symplectic geometry.  We first describe contact versions of the algebra of observables and Hamiltonian group actions, and give a construction of a Hilbert space of sections of a ``quantum bundle'' in the tradition of Kirillov-Kostant quantization.  Tools from CR geometry play a significant role in this construction; in particular, this approach applies to Sasakian manifolds.

The second approach is analogous to the use of Spin$^c$ (almost complex) quantization in symplectic geometry as a model for geometric quantization in the K\"ahler case \cite{GGK,SJ}: using a ``compatible'' almost CR structure, we construct an odd first order differential operator $\dirac$ that reduces, in the case of a strongly pseudoconvex CR manifold, to the operator $\dirac = \sqrt{2}(\dbbar+\dbbar^*)$, where $\dbbar$ is the tangential Cauchy-Riemann operator determined by the CR structure.  The operator $\dirac$ is not elliptic, but if a Lie group $G$ acts on $M$ transverse to the contact distribution, then $\dirac$ will be transversally elliptic, and we can give a formula for its index similar to the Riemann-Roch formula in the symplectic case.  

Let $(M,E)$ be a compact cooriented contact manifold.  A choice of contact form is given by a non-vanishing section $\theta$ of the annihilator line bundle $E^0\subset T^*M$.  (By assumption, $E^0$ is oriented, and hence, trivial.)  The subbundle $E=\ker\theta\subset TM$ is a contact distribution if and only if $\mu_\theta = \theta\wedge d\theta^n/n!$ defines a volume form on $M$.  If a compact Lie group $G$ acts on $M$ preserving $E$, the contact form $\theta$ can be assumed to be $G$-invariant by averaging, allowing us to define the contact momentum map $\Phi_\theta:M\to\mathfrak{g}^*$ given by
\[
\langle \Phi_\theta,X\rangle = \theta(X_M)
\]
for all $X\in\mathfrak{g}$, where $X_M$ is the vector field generated by the infinitesimal action of $X$ on $M$.  The contact form also determines a Jacobi structure on $M$ as follows: any vector field on $M$ is determined uniquely by its pairings with $\theta$ and $d\theta$; in particular, the Reeb vector field $\xi$ is defined by $\theta(\xi)=1$ and $\iota(\xi)d\theta = 0$.  This allows us to define a map $\Lambda^\#:T^*M\to E\subset TM$ by declaring that, for any $\eta\in T^*M$, we have
\[
\theta(\Lambda^\#\eta)=0 \quad\text{ and }\quad \iota(\Lambda^\#\eta)d\theta = \eta(\xi)\theta-\eta.
\]
Each $f\in C^\infty(M)$ is then associated to the Hamiltonian vector field $X_f = \Lambda^\#df + f\xi$, and the Jacobi bracket on $C^\infty(M)$ is given by $\{f,g\} = X_f\cdot g - g\xi\cdot f$.  For any $f\in C^\infty(M)$, the associated Hamiltonian vector field satisfies $\mathcal{L}(X_f)\theta = (\xi\cdot f)\theta$, so that $X_f$ is a contact vector field (see \cite{Lich1}).  We see that whenever $\xi\cdot f=0$, $X_f$ preserves the contact form, and hence the volume form $\mu_\theta$.
\begin{proposition}
The space $C_b^\infty(M) = \{f\in C^\infty(M)|\xi\cdot f = 0\}$ is a Lie subalgebra of $(C^\infty(M),\{\cdot,\cdot\})$, and the Jacobi bracket on $C^\infty(M)$ restricts to a Poisson bracket on $C^\infty_b(M)$.
\end{proposition}
In particular, since $\theta$ is preserved by the $G$-action, we can show that the momentum map components $\Phi^X_\theta = \langle \Phi_\theta, X\rangle \in C^\infty(M)$ satisfy $\xi\cdot \Phi^X_\theta = 0$ for all $X\in\mathfrak{g}$:
\begin{theorem}
Suppose a compact Lie group $G$ acts on a compact contact manifold $M$ preserving the contact form $\theta$.  With respect to the Jacobi structure determined by $\theta$, we have:
\begin{enumerate}
\item The map $\mathfrak{g}\to C^\infty(M)$ given by $X\mapsto \Phi^X_\theta$ is a Lie algebra homomorphism.
\item The Hamiltonian vector field associated to the function $\Phi^X_\theta$ is equal to $X_M$.
\end{enumerate}
\end{theorem}
In other words, the momentum map components span a Lie subalgebra of $C^\infty_b(M)\subset C^\infty(M)$, and the diagram of Lie algebra homomorphisms
\[
\xymatrix{\mathfrak{g}\ar[r]\ar[rd] &C^\infty(M)\ar[d]\\
 & \mathfrak{X}_{ham}(M)}
\]
commutes, where $\mathfrak{X}_{ham}(M)$ denotes the space of contact Hamiltonian vector fields on $M$.  

To define a quantization of the contact manifold $(M,\theta)$, we make use of the notion of a {\em quantum bundle} from \cite{DT}.  Since $\theta$ is a contact form, the 2-form $\Omega = -d\theta|_{E\otimes E}$ defines a symplectic structure on the subbundle $E=\ker\theta$.  A Hermitian line bundle with connection $\pi:(\mathbb{L},h,\nabla)\to (M,E,\Omega)$ is called a quantum bundle if the restriction of the curvature form of $\nabla$ to $E\otimes E$ is equal to $i\Omega$.  We can then construct the Hilbert space $\mathcal{H} = \Gamma_{L^2}(M,\mathbb{L})$ given by the $L^2$ completion of the space of smooth sections of $\mathbb{L}$ with respect to the inner product
\[
\langle s_1,s_2\rangle = \int_M h(s_1,s_2)\mu_\theta.
\]
It is then straightforward to check that the assignment
\[
f\mapsto \nabla_{X_f} + i\pi^*f
\]
defines a Lie algebra homomorphism from $C^\infty_b(M)$ to the space of Hermitian operators on $\mathcal{H}$.  In particular, we obtain a representation of the Lie algebra $\mathfrak{g}$ via the momentum map components $\Phi^X_\theta$.  Note however that the constant functions on $M$ do not correspond to multiples of the identity operator, since the contact Hamiltonian vector field associated to the constant $c$ is $X_c = c\xi$.

As in the symplectic case, we wish to reduce the size of the Hilbert space $\mathcal{H}$ by applying a polarization.  In this paper we consider ``compatible CR structures'' as the contact analogue of a complex polarization.  It then becomes natural to consider the case that $\mathbb{L}$ is a CR-holomorphic line bundle, and define our quantization to be the subspace of CR-holomorphic sections in $\mathcal{H}$.  It is possible to equip the trivial bundle $\mathbb{L}=M\times\C$ with the structure of a CR-holomorphic quantum bundle, so that the quantization of $M$ becomes the space of CR-holomorphic sections on $M$ (which in this case can be identified with the CR-holomorphic functions on $C^\infty(M)$.  This agrees with the ``fairly canonical'' quantization of Boutet de Monvel-Guillemin-Sternberg \cite{BMG,GS2} in the case where $M$ is the boundary of a strongly pseudoconvex complex domain.  (It may also be interesting to consider Legendrian foliations as an analogue of real polarizations, but we do not consider this problem here.)

The assumption of the existence of a compatible CR structure implies that $M$ is a strongly pseudoconvex CR manifold of hypersurface type, where the contact form $\theta$ determines a pseudo-Hermitian structure on $M$.  However, in general we can only expect a contact manifold to admit an almost CR structure.  As noted above, the contact distribution $E\subset TM$ is a symplectic subbundle, and we can choose a complex structure on the fibres of $E$ compatible with the symplectic structure.  The resulting splitting $E\otimes \C = E_{1,0}\oplus E_{0,1}$ determines an almost CR structure on $M$.  As in \cite{F}, using a compatible metric $g$ and connection $\nabla$ we can construct an odd first-order differential operator $\dirac$ acting on sections of $\S = \Lambda E_{0,1}^*$ similar to the Spin$^c$-Dirac operator associated to an almost Hermitian structure.  We use the metric to define a Clifford action of the bundle $\Cl(E)$ (defined by $\Cl(E)_x = \Cl(E^*_x,g|_{E^*_x})$) on $\S$, and define
\[
\dirac = \mathbf{c}\circ\pi_{E^*}\circ\nabla,
\]
where $\mathbf{c}$ denotes the Clifford multiplication, and $\pi_{E^*}:T^*M\to E^*$. We also can twist this construction by a quantum bundle $\mathbb{L}$ to obtain an operator $\dob_\mathbb{L}$ acting on sections of $\S\otimes\mathbb{L}$.
When the almost CR structure is integrable, so that $M$ is again a strongly pseudoconvex CR manifold, we can take our connection to be the Tanaka-Webster connection (see \cite{DT}).  In this case, we have
\begin{theorem}
On a strongly pseudoconvex CR manifold, if we define $\dirac$ using the Tanaka-Webster connection, then
\[
\dirac = \sqrt{2}\left(\dbbar+\dbbar^*\right),
\]
where $\dbbar:\mathcal{A}^{0,q}(M)\to \mathcal{A}^{0,q+1}(M)$ is the tangential Cauchy-Riemann operator.
\end{theorem}
When $\mathbb{L}$ is CR-holomorphic, we obtain an analogous result for the twisted operator $\dob_\mathbb{L}$, where $\dbbar$ is replaced in the above formula by the $\overline{\partial}_\mathbb{L}$ operator defining the CR-holomorphic structure on $\mathbb{L}$.  We use $\dirac$ to define the index-theoretic quantization $Q(M)$ given by  the $\Z_2$-graded space $Q(M)=\ker(\dob_\mathbb{L}^+)\oplus\ker(\dob_\mathbb{L}^-)$ (where $\dob_\mathbb{L}^\pm$ denotes the restrictions of $\dob_\mathbb{L}$ to even/odd forms).  

The operator $\dob_\mathbb{L}$ is not elliptic, since its principal symbol vanishes along the annihilator $E^0\subset T^*M$.  If a compact Lie group $G$ acts on $M$ preserving $\theta$, $g$, and $\nabla$, then $\dob_\mathbb{L}$ will commute with the $G$-action, and $Q(M)$ becomes a virtual $G$-representation, which in general is infinite-dimensional.  However, if the action of $G$ is transverse to the contact distribution, then $\dob_\mathbb{L}$ will be $G$-transversally elliptic, and the character of this representation, given by the equivariant index of $\dob_\mathbb{L}$, is defined as a generalized function (distribution) on $G$ \cite{AT}, and a cohomological formula for the index was given in \cite{F}.  Near the identity element in $G$, the index of $\dob_\mathbb{L}$ is given, for $X\in\mathfrak{g}$ sufficiently small, by
\[
 \ind^G(\dob_\mathbb{L})(e^X) = \frac{1}{(2\pi i)^n}\int_M \Td(E,X)\Ch(\mathbb{L},X)\theta\wedge\delta_0(d\theta - \theta(X_M)),
\]
with similar formulas near other elements of $G$.   Here, $\delta_0$ denotes the Dirac delta distribution on $\R$, so that $\theta\wedge\delta_0(d\theta-\theta(X_M))$ is an equivariant differential form with generalized coefficients; this form and its properties are explained in \cite{F}.

\section{Geometric quantization of symplectic manifolds}
Since the material presented here is quite standard, we will try to be brief, and refer the reader to the texts \cite{GGK,BW2,Wood} for details.  Let $(M,\omega)$ be a compact symplectic manifold. The  classical dynamics are given by Hamilton's equations: to any Hamiltonian function $H\in C^\infty(M)$ we can associate the unique vector field $X_H\in \mathfrak{X}(M)$ satisfying
\[dH = \iota(X_H)\omega.\]
Such vector fields are {\em symplectic}, in the sense that the flow of $X_H$ preserves the symplectic form $\omega$.  Moreover, the integral curves of $X_H$ lie in level sets of $H$. The algebra of observables is the Poisson algebra $C^\infty(M)$, equipped with the Poisson bracket $\{f,g\} = \omega(X_g,X_f)$.
We now suppose that a compact Lie group $G$ acts on $M$, preserving $\omega$.  This gives us a map
\begin{align*}
 \mathfrak{g}&\rightarrow \mathfrak{X}_{symp}(M)\\
  X&\mapsto X_M,
\end{align*}
where $\mathfrak{X}_{symp}(M)$ denotes the space of symplectic vector fields, and $X_M$ is the vector field generated by the infinitesimal action of $\mathfrak{g}$ on $M$.  The action of $G$ is called {\em Hamiltonian} if this map factors through the map $C^\infty(M)\rightarrow \mathfrak{X}_{symp}(M)$ given by associating a function to its Hamiltonian vector field.
\begin{definition}
 A {\bf momentum map} is an equivariant map $\Phi:M\rightarrow \mathfrak{g}^*$ such that
for each $X\in \mathfrak{X}(M)$, the pairing $\Phi^X = \langle\Phi,X\rangle$ satisfies
\begin{equation}\label{mmsymp}
 d\Phi^X = \iota(X_M)\omega.
\end{equation}
\end{definition}
Such a momentum map exists if and only if the action of $G$ on $(M,\omega)$ is Hamiltonian \cite{ACdS,GS1}).  The desired mapping $\mathfrak{g}\rightarrow C^\infty(M)$ is given by $X\mapsto \Phi^X$.  We note that this map is a Lie algebra homomorphism with respect to the Lie algebra structure on $C^\infty(M)$ given by the Poisson bracket.
The functions $\{\Phi^X|X\in \mathfrak{g}\}$ thus span a Lie subalgebra of $C^\infty(M)$. 

To our symplectic manifold $(M,\omega)$ we wish to associate a Hilbert space $\mathcal{H}$, such that the action of $G$ on $M$ corresponds to a representation of $G$ on $\mathcal{H}$.  Moreover, classical `observables' should correspond to quantum ones: there should be an algebra of skew-Hermitian operators $A_X$ on $\mathcal{H}$ and Lie algebra homomorphism to this algebra (with respect to the commutator bracket) from the algebra generated by the momentum map components $\Phi^X$ (with respect to the Poisson bracket).

Suppose we are given a Hamiltonian $G$-space $(M,\omega,\Phi)$, such that the equivariant cohomology class of $\omega(X) = \omega - \Phi(X)$ is integral.  Then there exists a $G$-equivariant complex line bundle $\pi:\mathbb{L}\rightarrow M$, equipped with $G$-invariant Hermitian metric $h$ and connection $\nabla$ with equivariant curvature form $F_\nabla(X) = i\omega(X)$. Such a line bundle $\mathbb{L}$ is known as a {\bf $G$-equivariant prequantum line bundle} for $(M,\omega,\Phi)$.
The action of $G$ on $M$ induces a linear action of $G$ on the space of sections of $\mathbb{L}$ by bundle automorphisms, and we obtain a unitary representation of $G$ on the Hilbert space 
\[
 \mathcal{H} = \Gamma_{L^2}(M,\mathbb{L})
\]
of $L^2$ sections of $\mathbb{L}$, with respect to the inner product
\begin{equation*}
 \langle s_1,s_2\rangle = \int_M h(s_1,s_2) \frac{\omega^n}{n!}.
\end{equation*}
From the infinitesimal action of $\mathfrak{g}$ on the space of sections, we obtain the desired correspondence $\Phi^X\mapsto A_X$ between classical and quantum observables via
\[A_X = \nabla_{X_M} + i\pi^*\Phi^X.\]

The Hilbert space we obtain in this way turns out to be too big (for example, in the non-compact case $M=T^*X$, we obtain $L^2(T^*X)$ rather than $L^2(X)$, as one would expect in the canonical Schr\"odinger quantization).  The standard way of cutting down the space of sections is to apply a {\em polarization}.  We will restrict ourselves to the case of a complex polarization, which is defined to be an integrable maximal isotropic subbundle $\mathcal{P}$ of  $TM\otimes\C$ such that $\mathcal{P}\cap \overline{\mathcal{P}}=0$.

In other words, a polarization is given by a complex structure on $M$ that is compatible with the symplectic structure.  The existence of a complex polarization is thus equivalent to having a K\"ahler structure on $M$.
A polarization determines a subspace of the space of $L^2$ sections of $\mathbb{L}$ by requiring $\nabla_{\overline{X}} s = 0$ for all $X\in\mathcal{P}$; these are the so-called {\em polarized} sections. The space of polarized sections is then a candidate for the space $Q(M)$. By \cite[Proposition 6.30]{GGK}, there is a unique holomorphic structure on $\mathbb{L}$ such that the (local) polarized sections of $\mathbb{L}$ are the (local) holomorphic sections of $\mathbb{L}$.  That is, the connection $\nabla$ preserves the metric, and satisfies $\nabla^{0,1} = \overline{\partial}_\mathbb{L}$, where $\nabla^{0,1} = \nabla|_{\overline{\mathcal{P}}}$.  The resulting quantization $Q(M)$ in this case is then given by the space of holomorphic sections of $\mathbb{L}$.
%\begin{remark}
%We can see here why it is advantageous to work with sections of a complex line bundle, rather than complex functions on $M$: the space $Q(M)=\Gamma_{hol}(M,\mathbb{L})$ is non-trivial if  $\mathbb{L}$ is sufficiently positive, while this is not the case for the space of holomorphic functions on the compact manifold $M$.
%Moreover, we can generalize slightly by replacing the space of holomorphic sections of $\mathbb{L}$ by the sheaf cohomology groups $H^k(M;\mathcal{O}(\mathbb{L}))$.  The definition of $Q(M)$ is then taken to be the $\mathbb{Z}_2$-graded vector space
%\begin{equation}
% Q(M) = \sum (-1)^kH^k(M;\mathcal{O}(\mathbb{L})),
%\end{equation}
%where the $(-1)^k$ term is being used as shorthand to indicate the grading. The corresponding action of $G$ on these groups is well-understood \cite{Bott}, and moreover, the groups $H^k(M;\mathcal{O}(\mathbb{L}))$ are isomorphic to the Dolbeault cohomology groups $H^{0,k}(M;\mathbb{L})$ of the complex of differential forms on $M$ with values in $\mathbb{L}$.  This observation will later allow us to think of $Q(M)$ in terms of equivariant index theorem.
%\end{remark}

\section{Geometric quantization of contact manifolds}
\subsection{Contact momentum maps}
Let $(M,E)$ be a compact contact manifold of dimension $2n+1$.  We will assume that the contact distribution $E$ is cooriented, so that there exists a global contact form $\theta\in\Gamma(M,E^0\setminus 0)$.  The contact form $\theta$ determines a splitting $T^*M = E^*\oplus E^0$ of the cotangent bundle, a trivialization $E^0 = M\times \R$, and an orientation on $M$ given by the volume form $\mu_\theta = \theta\wedge d\theta^n/n!$.

We suppose that a compact Lie group $G$ acts on $M$ by contactomorphisms; by averaging, we may assume that the contact form $\theta$ is $G$-invariant.
\begin{definition}
The {\bf contact momentum map} associated to the contact form $\theta$ is the map $\Phi_\theta:M\rightarrow \mathfrak{g}^*$ such that for any $X\in \mathfrak{g}$, we have
\begin{equation}\label{mmcont}
 \langle\Phi_\theta,X\rangle = \theta(X_M).
\end{equation}
\end{definition}
\begin{remark}
 The contact momentum map defined above does of course depend on the choice of contact form $\theta$.  For further discussion of the properties of contact momentum maps, see \cite{ler}.
\end{remark}
We note that the momentum map components $\Phi^X_\theta = \theta(X_M)$ satisfy similar properties to the components of a symplectic momentum map.  In particular, $\Phi^X_\theta$ is the `Hamiltonian' function associated to the vector field $X_M$, in the sense that, by the invariance of $\theta$, we have
\[
 d\Phi^X_\theta = d\iota(X_M)\theta = -\iota(X_M)d\theta = \iota(X_M)\Omega,
\]
where $\Omega = -d\theta$ restricts to a symplectic structure on the fibres of the contact distribution $E$ determined by $\theta$.
\subsection{The Jacobi algebra}
Given an action of $G$ on $(M,\theta)$ leaving $\theta$ invariant, the vector fields $X_M$ generated by the Lie algebra elements $X\in\mathfrak{g}$ are contact, since $\mathcal{L}(X_M)\theta = 0$ for all $X\in\mathfrak{g}$. (In general, a contact vector field $V$ satisfies $\mathcal{L}(V)\theta = f\theta$ for some $f\in C^\infty(M)$; the vector fields that preserve the contact form are characterized physically in \cite{Raj} as the {\em incompressible} vector fields, since they also preserve the volume form $\mu_\theta$.)  As with symplectic geometry, there is a standard notion of a Hamiltonian vector field associated to each function on a contact manifold: given $f\in C^\infty(M)$, the contact Hamiltonian vector field $X_f$ is the unique vector field such that $\theta(X_f) = f$ and $\iota(X_f)d\theta = (\xi\cdot f)\theta - df$.  Moreover, a choice of contact form determines a Lie algebra structure on $C^\infty(M)$ via the {\em Jacobi bracket}.   Jacobi structures were first developed (independently) by Kirillov \cite{Kir2} and Lichnerowicz \cite{Lich}; our primary reference for this section is the article \cite{Marle}, although for the definition of contact Hamiltonian vector fields and the resulting Jacobi bracket on functions (without reference to general Jacobi structures), see \cite{Lich1}.
\begin{remark}
A bracket on smooth functions called the {\em Lagrange bracket} is defined (in terms of coordinates) in \cite{Raj}.  It is straightforward to check that, up to a sign convention, the Jacobi bracket defined below reduces to the Lagrange bracket in a Darboux chart.  Thus, the ``generalized Poisson algebra'' defined in\cite{Raj} is simply the usual Jacobi algebra structure for the standard contact structure on $\R^{2n+1}$.
\end{remark}
\begin{definition}
 A {\bf Jacobi structure} on a manifold $M$ is a bracket $\{\cdot,\cdot\}$ on $C^\infty(M)$ that is skew-symmetric, satisfies the Jacobi identity, and is local, in the sense that the support of $\{f,g\}$ is contained in the intersection of the supports of $f$ and $g$.
\end{definition}
A Jacobi structure is equivalent to the existence of a bivector field $\Lambda\in\Gamma(M,\Lambda^2(TM))$ and a vector field $\xi\in\mathfrak{X}(M)$ such that
\begin{equation*}
%\label{jacobi}
 [\xi,\Lambda] = \mathcal{L}(\xi)\Lambda = 0 \quad\text{and}\quad [\Lambda,\Lambda] = 2\xi\wedge\Lambda,
\end{equation*}
where $[\cdot,\cdot]$ denotes the Schouten bracket.  The relationship between the Jacobi bracket and the data $(\Lambda,\xi)$ is given by
\begin{equation*}
%\label{bracket}
 \{f,g\} = \Lambda(df,dg)+\iota(\xi)(f\,dg-g\,df).
\end{equation*}
Given a contact manifold $(M,E)$ equipped with contact form $\theta$, the vector field $\xi$ is given by the {\em Reeb field}, defined to be the unique vector field such that
\[
 \iota(\xi)\theta = 1 \quad\text{and}\quad\iota(\xi)\Omega = 0,
\]
where $\Omega = -d\theta$.
The contact form also determines a map $\Lambda^{\#}:T^*M\rightarrow TM$ such that, for any $\eta\in T^*M$, we have
\begin{equation*}
%\label{lammap}
 \theta(\Lambda^\#(\eta)) = 0\quad\text{and}\quad\iota(\Lambda^\#(\eta))\Omega = \eta - (\eta(\xi))\theta.
\end{equation*}
We note that the image of the map $\Lambda^\#$ is contained in the contact distribution $E$, by the first of the above two conditions.  Finally, we can define $\Lambda\in\Gamma(M,\Lambda^2(TM))$ by
\begin{equation*}
%\label{lambda}
 \Lambda(\eta,\zeta) = \iota(\Lambda^\#(\eta))\zeta = -\iota(\Lambda^\#(\zeta))\eta.
\end{equation*}
From the Jacobi structure associated to the contact form $\theta$, we obtain a Lie algebra structure on $C^\infty(M)$, as well as a notion of Hamiltonian vector field:
\begin{definition}
 For any $f\in C^\infty(M)$, the {\bf Hamiltonian vector field} associated to $f$ is the vector field
\begin{equation}\label{ham}
 X_f = \Lambda^\#(df) + f\xi.
\end{equation}
\end{definition}
The following facts can be found in \cite{Lich1}, although the first fact is true for Jacobi structures in general:
\begin{proposition}\label{hamprop}
For any $f\in C^\infty(M)$, the associated Hamiltonian vector field $X_f$ satisfies the following properties:
\begin{enumerate}
\item  The map $f\mapsto X_f$ is a Lie algebra homomorphism: for any $f,g\in C^\infty(M)$, we have
\[
 X_{\{f,g\}} = [X_f,X_g].
\]
\item For any $g\in C^\infty(M)$, $X_f\cdot g = \{f,g\}+(\xi\cdot f)g$.
\item $\iota(X_f)\Omega = df-(\xi\cdot f)\theta$.
\item $X_f$ is a contact vector field: $\mathcal{L}(X_f)\theta = (\xi\cdot f)\theta$.
\end{enumerate}
\end{proposition}
%\begin{proof}
%The first property is well-known and appears for example in \cite{Marle}.
%The proof of each of the remaining properties is obtained by a straight-forward calculation.  For the second, we have
%\[
%X_f\cdot g = \iota(\Lambda^\#(df)+f\xi)dg = \Lambda(df,dg) + f\xi\cdot g = \{f,g\}+g\xi\cdot f.
%\]
%For the third, we have
%\[
%\iota(X_f)\Omega = \iota(\Lambda^\#(df))\Omega + f\iota(\xi)\Omega = df - (\iota(\xi)df)\theta,
%\]
%and the fourth now follows from the third by Cartan's formula, since $\iota(X_f)d\theta = -\iota(X_f)\Omega$, and $d\iota(X_f)\theta = df$.
%\end{proof}
\subsection{The Poisson algebra}
From the above proposition, we see that the image of the homomorphism $C^\infty(M)\rightarrow\mathfrak{X}(M)$ given by \eqref{ham} is contained in the Lie subalgebra of contact vector fields.  (This was essentially the goal of the construction given in \cite{Lich1}.)  Moreover, we note that in each case, the failure of $(C^\infty(M),\{\cdot ,\cdot\})$ to behave like a Poisson algebra is indicated by the presence of the term $\xi\cdot f$.  We therefore might ask what can be said about those functions for which $\xi\cdot f = 0$.
For any manifold $M$ equipped with a closed two-form $\Omega$, we have the associated Poisson algebra \cite{GGK}
\[
\mathcal{P}(M,\Omega) = \{(f,X)\in C^\infty(M)\times\mathfrak{X}(M)|df = \iota(X)\Omega\}.
\]
The bracket is given by $[(f,X),(g,Y)] = (\frac{1}{2}(Y\cdot f-X\cdot g),[X,Y])$, and the Poisson algebra acts on $M$ via $(f,X)\mapsto X$.
Of course, if $\Omega$ is symplectic, then $\mathcal{P}(M,\Omega)$ is isomorphic to $C^\infty(M)$,
since each $f$ is associated to a unique Hamiltonian vector field $X_f$.

Let us suppose instead that $(M,\theta)$ is a contact manifold, and consider the Poisson algebra $\mathcal{P}(M,\Omega)$ with respect to the two-form $\Omega = -d\theta$.  For any $f\in C^\infty(M)$, we can consider the pair $(f,X_f)$, where $X_f = \Lambda^\#(df) + f\xi$, as above.  By Proposition \ref{hamprop}, we see that $\iota(X_f)\Omega = df$ if and only if $\xi\cdot f = 0$.
\begin{lemma}\label{xidot}
For any $f\in C^\infty(M)$, $[\xi,X_f] = X_{\xi\cdot f}$.
\end{lemma}
\begin{proof}
We simply check that $\iota([\xi,X_f])\theta = \xi\cdot f$ and $\iota([\xi,X_f])\Omega = d(\xi\cdot f) - \xi\cdot(\xi\cdot f)\theta$.
\end{proof}
\begin{lemma}\label{xijac}
For any $f,g\in C^\infty(M)$, $\xi\cdot\{f,g\} = \{\xi\cdot f, g\} +\{f,\xi\cdot g\}$.
\end{lemma}
\begin{proof}
Using Lemma \ref{xidot}, we see that
\begin{align*}
\xi\cdot\{f,g\} & = \xi\cdot (X_f\cdot g) - \xi\cdot ((\xi\cdot f)g)\\
& = X_f\cdot\xi\cdot g + X_{\xi\cdot f} g - g\xi\cdot(\xi\cdot f) -(\xi\cdot f)(\xi\cdot g)\\
& = \{f,\xi\cdot g\} + \{\xi\cdot f,g\}.
\end{align*}
\end{proof}
\begin{remark}
Using the above two lemmas, it is straightforward to verify directly that $[X_f,X_g]=X_{\{f,g\}}$ by computing the contractions of each with $\theta$ and $\Omega$.
\end{remark}
\begin{definition}
We denote by $\mathcal{P}_b(M,\Omega)$ the subset of $\mathcal{P}(M,\Omega)$ given by
\[
\mathcal{P}_b(M,\Omega) = \{(f,X_f)\in C^\infty(M)\times\mathfrak{X}(M)|\xi\cdot f = 0\}.
\]
\end{definition}
We note that $\mathcal{P}_b(M,\Omega)$ is a proper subset of $\mathcal{P}(M,\Omega)$, since $(f,X_f+g\xi)\in\mathcal{P}(M,\Omega)$ for any $g\in C^\infty(M)$.  Moreover, we have the following:
\begin{proposition}
The set $\mathcal{P}_b(M,\Omega)$ is a Lie subalgebra of $\mathcal{P}(M,\Omega)$.
\end{proposition}
\begin{proof}
For any $f,g\in C^\infty(M)$, we see using Proposition \ref{hamprop} that
\[
[(f,X_f),(g,X_g)] = (\{f,g\}+\frac{1}{2}(g\xi\cdot f - f\xi\cdot g),[X_f,X_g]).
\]
Since $[X_f,X_g] = X_{\{f,g\}}$, whenever $\xi\cdot f = \xi\cdot g = 0$ we have $[(f,X_f),(g,X_g)] = (\{f,g\},X_{\{f,g\}})$  and by Lemma \ref{xijac}, $\xi\cdot \{f,g\}=0$, so that the pair $(\{f,g\},X_{\{f,g\}})$ belongs to $\mathcal{P}_b(M,\Omega)$.
\end{proof}
Using the above, and the fact that $\xi\cdot f = 0$ if and only if $\iota(X_f)\Omega = df$, we obtain the following:
\begin{proposition}
If we define the subsets $C^\infty_b(M) = \{f\in C^\infty(M)|\xi\cdot f = 0\}$ and $\mathfrak{X}_b(M) = \{X_f\in\mathfrak{X}_{ham}(M)|\iota(X_f)\Omega = df\}$, then we have:
\begin{enumerate}
\item The space $C^\infty_b(M)$ is a Lie subalgebra of $(C^\infty(M),\{\cdot , \cdot\})$.
\item The space $\mathfrak{X}_b(M)$ is a Lie subalgebra of $(\mathfrak{X}(M),[\cdot ,\cdot])$.
\item We have Lie algebra isomorphisms $\mathcal{P}_b(M,\Omega)\cong C^\infty_b(M) \cong \mathfrak{X}_b(M)$.
\end{enumerate}
\end{proposition}

In particular, the above tells us that the Jacobi subalgebra $C^\infty_b(M)\subset C^\infty(M)$ is in fact a Poisson algebra.
 Let us further denote by $\mathfrak{X}_{symm}(M,\theta) = \{X\in\mathfrak{X}(M)|\mathcal{L}(X)\theta = [X,\xi]=0\}$ the Lie algebra of infinitesimal symmetries of $(M,\theta)$.  Following \cite{GGK}, we have:
\begin{proposition}
The map $\mathcal{P}_b(M,\Omega)\rightarrow \mathfrak{X}_{symm}(M,\theta)$ given by $(f,X_f)\mapsto X_f$ is an isomorphism of Lie algebras.
\end{proposition}
\begin{proof}
By Proposition \ref{hamprop}, we see that the action of the pair $(f,X_f)\in\mathcal{P}(M,\Omega)$ on $M$ preserves the contact form, since $\mathcal{L}(X_f)\theta = (\xi\cdot f)\theta = 0$.  Moreover, we have
$[X_f,\xi] = 0$ by Lemma \ref{xidot},
so that $X_f$ is an infinitesimal symmetry of $(M,\theta)$.  Conversely, choose any $X\in \mathfrak{X}_{symm}(M,\theta)$.  Let us write $X=Y+f\xi$, where $f=\theta(X)$ and $Y=X-\theta(X)\xi\in E=\ker{\theta}$. Since $\mathcal{L}(X)\theta = 0$, we have 
\[
0 = \mathcal{L}(Y+f\xi)\theta = \iota(Y+f\xi)d\theta + d\iota(Y+f\xi)\theta = \iota(Y)d\theta + df,
\]
and therefore, $df = -\iota(Y)d\theta = \iota(Y)\Omega$.  Since $\iota(Y)\Omega = df$ and $\iota(Y)\theta = 0$, it follows that $Y=\Lambda^\#(df)$, and thus, $X = X_f$.
\end{proof}

Let us now consider the case where $(M,\theta)$ is a Boothby-Wang fibration \cite{BW}.  That is, $(M,\theta)$ is a principal $U(1)$-bundle over a symplectic manifold $(B,\omega)$, with connection 1-form $\theta$ (identifying $\mathfrak{u}(1)$ with $\R$).   The symplectic manifold $(B,\omega)$ is then  prequantizable, and $\mathbb{L}=M\times_{U(1)}\C$ is the associated prequantum line bundle.  Given an action of $G$ on $(M,\theta)$ preserving $\theta$, the prequantization condition becomes $\pi^*\omega = -d\theta$ and $\pi^*\Phi^X = \iota(X_M)\theta$.

As outlined in \cite{GGK}, the traditional Kirillov-Kostant approach is to start with a Hamiltonian action of $G$ on $(B,\omega)$, and try to lift the infinitesimal action to $M$ such that $\pi^*\Phi^X = \theta(X_M)$.  However, one can in fact lift the action of the entire Poisson algebra $C^\infty(B)$ to $M$: given $f\in C^\infty(B)$, let $X_f$ be its associated (symplectic) Hamiltonian vector field.  The action of $f$ on $M$ is then given by
\begin{equation}\label{lift}
f\mapsto X_f^{hor} + \pi^*f\cdot \xi,
\end{equation}
where $X_f^{hor}$ denotes the horizontal lift of $X_f$ with respect to the connection $\theta$, and the Reeb field $\xi$ is the infinitesimal generator of the $U(1)$ action.  By \cite[Proposition 6.17]{GGK}, the Poisson algebra of $(B,\omega)$ is isomorphic via the above map to $\mathfrak{X}_{symm}(M)$, and thus, $\mathcal{P}(B,\omega)$ is isomorphic to $\mathcal{P}_b(M,\Omega)$.  Moreover, we see that the vector field \eqref{lift} is the Hamiltonian vector field (in the Jacobi sense, given by \eqref{ham}) associated to $\pi^*f$.
\begin{remark}
One advantage of our approach is that the algebra $\mathcal{P}_b(M,\Omega)$ makes sense even when $(M,\theta)$ is not a regular contact manifold (that is, when the Reeb field corresponding to $\theta$ does not generate a free circle action), and therefore can be applied in settings where no regular contact structure exists.  Moreover, we notice that when the lifts of the momentum map components $\Phi^X$ satisfy the prequantization condition, they exactly coincide with the components of the contact momentum map.

In Section \ref{crpolar} below, we will see that the trivial line bundle $\mathbb{L}=M\times\C$ serves as a contact version of the prequantum line bundle.  When $(M,\theta)$ is a prequantum circle bundle, we note that $\mathbb{L}/U(1) = M\times_{U(1)}\C$ is a prequantum line bundle for the symplectic manifold $M/U(1)$.
\end{remark}
\subsection{Contact momentum maps revisited}
For any $X\in \mathfrak{g}$ we have the function $\Phi^X_\theta\in C^\infty(M)$ given in terms of the contact momentum map.  Continuing the analogy with symplectic geometry, we have the following:
\begin{theorem}
 Suppose a compact Lie group $G$ acts on a compact contact manifold preserving a chosen contact form $\theta$, and let $\Phi_\theta:M\rightarrow \mathfrak{g}^*$ denote the corresponding contact momentum map.  With respect to the Jacobi structure defined by $\theta$, we have the following:
\begin{enumerate}
 \item The map $\mathfrak{g}\rightarrow C^\infty(M)$ given by $X\rightarrow \Phi^X_\theta$ is a Lie algebra homomorphism.
 \item The Hamiltonian vector field associated to $\Phi^X_\theta$ is equal to $X_M$.
\end{enumerate}
\end{theorem}
\begin{proof}
 Let $\Omega = -d\theta$, and note that $\iota(\xi)\Omega = 0$, and $\iota(X_M)\Omega = d\Phi^X_\theta$.  The Jacobi bracket is given by
\begin{align*}
 \{\Phi^X_\theta,\Phi^Y_\theta\} & = \Lambda(d\Phi^X_\theta,d\Phi^Y_\theta) + \iota(\xi)(\Phi^X_\theta \, d\Phi^Y_\theta - \Phi^Y_\theta\, d\Phi^X_\theta)\\
& = \Lambda(\iota(X_M)\Omega,\iota(Y_M)\Omega) + \Phi^X_\theta\Omega(Y_M,\xi) - \Phi^Y_\theta\Omega(X_M,\xi)\\
& = \Omega(Y_M,\Lambda^\#(\iota(X_M)\Omega))\\
& = -\iota(Y_M)[\iota(\Lambda^\#(\iota(X_M)\Omega))\Omega]\\
& = -\iota(Y_M)(\iota(X_M)\Omega - \Omega(X_M,\xi),\theta)\\
& = \Omega(Y_M,X_M),
\end{align*}
while the component of $\Phi_\theta$ in the direction of $[X,Y]$ is given by
\begin{align*}
 \Phi^{[X,Y]}_\theta & = \iota([X,Y]_M)\theta  = \iota([X_M,Y_M])\theta\\
 & = [\mathcal{L}(X_M),\iota(Y_M)]\theta\\
 & = \mathcal{L}(X_M)(\iota(Y_M)\theta) + \iota(Y_M)(\mathcal{L}(X_M)\theta)\\
 & = \iota(X_M)d(\iota(Y_M)\theta)\\
 & = -\iota(X_M)\iota(Y_M)d\theta = \Omega(Y_M,X_M),
\end{align*}
using the invariance of $\theta$.  This establishes the first point.  For the second, we note that the Hamiltonian vector field associated to $f = \Phi^X_\theta$ is given by
\begin{align*}
 X_f &= \Lambda^\#(d\Phi^X_\theta) + \Phi^X_\theta \xi\\
 & = \Lambda^\#(\iota(X_M)\Omega) + \Phi^X_\theta \xi.
\end{align*}
We now compute $\iota(X_f)\theta$ and $\iota(X_f)\Omega$.  We have
\[
 \iota(X_f)\theta = \iota(\Lambda^\#(\iota(X_M)\Omega))\theta + \iota(\Phi^X_\theta \xi)\theta = \Phi^X_\theta = \iota(X_M)\theta,
\]
and
\[
 \iota(X_f)\Omega = \iota(\Lambda^\#(\iota(X_M)\Omega))\Omega + \iota(\Phi^X_\theta\xi)\Omega =\iota(X_M)\Omega-\Omega(X_M,\xi)\theta = \iota(X_M)\Omega.
\]
\end{proof}
Thus, we see that any group action preserving the contact distribution is Hamiltonian, in the sense that, once an invariant contact form has been chosen, the map $\mathfrak{g}\rightarrow \mathfrak{X}_{cont}(M)$ factors through $C^\infty(M)$, and the image of this map is contained in the set of contact Hamiltonians.  As noted above, for an arbitrary $f\in C^\infty(M)$ the associated Hamiltonian vector field satisfies 
\begin{equation}\label{contvf}
\mathcal{L}(X_f)\theta = (\xi\cdot f)\theta.
\end{equation}  Since $\theta$ is $G$-invariant, for any $X\in \mathfrak{g}$ we have $\mathcal{L}(X_M)\theta = 0$.  Since $X_M$ is the Hamiltonian vector field associated to $\Phi^X_\theta$, we may deduce that
\[
 \xi\cdot \Phi^X_\theta = 0,
\]
so that $(\Phi_\theta^X,X_M)\in\mathcal{P}_b(M,\Omega)$.  Hence, we can consider the quantization of the contact manifold $(M,\theta)$ equipped with a group of symmetries $G$, in terms of the smaller Lie subalgebra spanned by the momentum map components $\Phi_\theta^X$.
\begin{remark}
As noted in the previous section, an advantage of working with the algebra $\mathcal{P}_b(M,\Omega)$ is that it makes sense when $M$ is not a Boothby-Wang fibration, although in this case this algebra is usually smaller.  When the dynamics of the Reeb field $\xi$ are more complicated, it is not immediately clear that it is even possible to find globally defined (non-constant) functions satisfying $\xi\cdot f = 0$.  However, the above results regarding the moment map components tells us that these functions always belong to $\mathcal{P}_b(M,\Omega)$.  This suggests that a contact manifold $M$ only admits a non-trivial action of a compact Lie group $G$ by contactomorphisms if there exist global non-constant solutions to the equation $\xi\cdot f = 0$, where $\xi$ is the Reeb field of an invariant contact form.

To see that it is possible to have a non-trivial group action preserving a contact form with a somewhat badly behaved Reeb field, we consider the following example of a contact structure on the 3-torus $\mathbb{T}^3$ due to Blair \cite{Blair} (who also shows that no regular contact structure on $\mathbb{T}^3$ exists):

On $\R^3$, we define the 1-form $\eta = \cos x_3 dx_1 + \sin x_3 dx_2$, which is a contact form invariant under the action $x_i\mapsto x_i+2\pi$, and thus it descends to a contact form $\theta = \cos\phi_3 d\phi_1 + \sin\phi_3 d\phi_2$ on $M=\mathbb{T}^3$.  The corresponding Reeb vector field is given by
\[
\xi = \cos\phi_3\frac{\partial}{\partial \phi_1} + \sin\phi_3\frac{\partial}{\partial \phi_2}.
\]
As noted by Blair, the integral curve of $\xi$ through $(0,0,\pi/3)$ is given by $t\mapsto (\frac{1}{2}t,\frac{\sqrt{3}}{2}t,\frac{\pi}{3})$, which is an irrational flow on the sub-2-torus $\phi_3=\pi/3$.  $M$ is therefore not a regular contact manifold.

On the other hand, we note that the contact distribution is spanned by the vector fields $X=\sin\phi_3\dfrac{\partial}{\partial \phi_1}-\cos\phi_3\dfrac{\partial}{\partial \phi_2}$ and $Y=\dfrac{\partial}{\partial \phi_3}$, and the action of $\mathbb{T}^2$ on $M$ given by
\[
(\alpha_1,\alpha_2)\cdot (\phi_1,\phi_2,\phi_3) = (\phi_1+\alpha_1,\phi_2+\alpha_2,\phi_3)
\]
preserves the contact distribution (and the contact form $\theta$).  Moreover, the action is transverse to the contact distribution, which will be relevant later when we consider an index-theoretic approach to contact quantization.  Note that in this example we have 
\[
C^\infty_b(M) = \left\{f\in\C^\infty(M)|\cos\phi_3\frac{\partial f}{\partial \phi_1}+\sin\phi_3\frac{\partial f}{\partial\phi_2}=0\right\}, 
\]
so that in particular any function of $\phi_3$ only belongs to $C^\infty_b(M)$.  With respect to the standard basis for $\R^2 \cong \Lie(\mathbb{T}^2)$ we have the moment map components $\Phi_\theta^{(1,0)} = \cos\phi_3$ and $\Phi_\theta^{(0,1)} = -\sin\phi_3$.
\end{remark}
\subsection{Quantum bundles}
Having established contact analogues of the symplectic description of a classical system, we now consider the corresponding construction of a quantum system.  We begin with {\em quantum bundles}, the generalization of prequantum line bundles described in \cite{DT}.

Let $M$ be a compact manifold, let $E\subset TM$ be a subbundle equipped with a symplectic form $\Omega$, and let $\pi:\mathbb{L}\rightarrow M$ be a Hermitian line bundle equipped with metric $h$ and connection $\nabla$.  

%To the connection $\nabla$ and any local section $s\in\Gamma(U\subset M,\mathbb{L}\setminus 0)$ we can associate a locally defined 1-form $\alpha(s)$ given by
%\[
% \iota(X)\alpha(s) = \frac{1}{i}\frac{\nabla_X s}{s}.
%\]
%There is then a unique connection form $\alpha\in\mathcal{A}^1(M,\mathbb{L}\setminus 0)$ such that for any local section $s\in\Gamma(U,\mathbb{L}\setminus 0)$, we have $s^*\alpha = \alpha(x)$.  The meaning of {\em connection form}, as defined in \cite{DT}, is the following: $\alpha$ should be $\operatorname{GL}(1,\C)$-invariant, and for any $x\in M$ and any $\operatorname{GL}(1,\C)$-equivariant map $\phi:\C\setminus 0\rightarrow \mathbb{L}_x\setminus 0$, we have
%\[
% \lambda^*j_x^*\alpha = \frac{1}{iz}dz,
%\]
%where $j_x:\mathbb{L}_x\setminus 0\hookrightarrow \mathbb{L}$ denotes inclusion.  There is then a unique closed 2-form $F_\nabla\in\mathcal{A}^2(M)$, the {\em curvature form} of $\nabla$, such that $\pi_0^*F_\nabla = d\alpha$, where $\pi_0$ denotes the restriction of $\pi:\mathbb{L}\rightarrow M$ to $\mathbb{L}\setminus 0$.

\begin{definition}
 We say that $(\mathbb{L},h,\nabla)\rightarrow (M,E,\Omega)$ is a {\bf quantum bundle} if the restriction of the curvature form of $\nabla$ to $E\otimes E$ is equal to $i\Omega$.
\end{definition}
Given a compact contact manifold $(M,E)$ and a choice of contact form $\theta$, we have the symplectic structure given by $\Omega = -d\theta|_{E\otimes E}$ on $E$.  Let us suppose that $(\mathbb{L},h, \nabla)$ is a quantum bundle over $(M,E,\Omega)$.  Then, following the symplectic case, we can consider the Hilbert space
\[
 \mathcal{H} = \Gamma_{L^2}(M,\mathbb{L}),
\]
with respect to the inner product
\[
 \langle s_1,s_2\rangle = \int_M h(s_1,s_2)\mu_\theta.
\]
We can define a map from $C^\infty(M)$ to the space of skew-Hermitian operators on $\mathcal{H}$ via the assignment
\begin{equation}\label{ftoA}
f\mapsto A_f = \nabla_{X_f}+i\pi^*f,
\end{equation}
where $\pi:\mathbb{L}\to M$.  In general this is not a Lie algebra homomorphism; however, we recall that the Jacobi bracket on $C^\infty(M)$ restricts to a Poisson bracket on $C^\infty_b(M) = \{f\in C^\infty(M)|\xi\cdot f = 0\}$.  Moreover, we have
\begin{proposition}
The restriction of the map \eqref{ftoA} to $C^\infty_b(M)$ is a Lie algebra homomorphism.
\end{proposition}
\begin{proof}
For any $f,g\in C^\infty_b(M)$, we have
\begin{align*}
[A_f,A_g] &= [\nabla_{X_f},\nabla_{X_g}] + i\pi^*\left(X_f\cdot g - X_g\cdot f\right)\\
&= \nabla_{[X_f,X_g]} + i\Omega(X_f,X_g)+i\pi^*\left(X_f\cdot g - X_g\cdot f\right)\\
&= \nabla_{X_{\{f,g\}}}+i\pi^*\{f,g\} = A_{\{f,g\}}.
\end{align*}
\end{proof}
In particular, we recall that the components of the contact momentum map $\Phi_\theta:M\rightarrow \mathfrak{g}^*$, belong to $C^\infty_b(M)$, so that for each $X\in\mathfrak{g}$, we can define the operator
\[
 A_X = \nabla_{X_M} + i\pi^*\Phi_\theta^X
\]
 on $\mathcal{H}$.  We note however that in the contact case, we cannot satisfy the ``normalization'' Dirac axiom, which requires that constant functions on $M$ correspond to multiples of the identity operator on $\mathcal{H}$: for each constant function $c\in C^\infty_b(M)$, we have $X_c = c\xi$.  This is reasonable from the point of view that the quantization of $M$ corresponds to the homogeneous quantization of the symplectization of $M$ (as in \cite{GS2}, for example), since constant functions on $M$ do not correspond to constant functions on the symplectization.  This is also consistent with the results obtained in \cite{Raj} using deformation quantization.

\subsection{CR polarizations}\label{crpolar}
As in the symplectic case, it is desirable to cut down the Hilbert space $\mathcal{H}$ to a smaller subspace.  Since our contact manifold $(M,E)$ is odd-dimensional, we cannot define a complex polarization on $M$.  Instead, we make the following definition:
\begin{definition}
 A subbundle $\mathcal{P}\subset T_\C M$ will be called a {\bf CR polarization} of the contact manifold $(M,E)$ provided that $\mathcal{P}$ is isotropic, formally integrable, $\mathcal{P}\cap \overline{\mathcal{P}} = 0$, and $\mathcal{P}\oplus\overline{\mathcal{P}} = E\otimes\C$.
\end{definition}
In other words, a CR polarization is simply a CR structure on $M$ whose Levi distribution is the contact distribution $E$.  Given such a polarization, we define the ``CR quantization'' of $(M, E)$ to be 
\[
Q(M) = \{s\in\Gamma_{L^2}(M,\mathbb{L})|\nabla_{\overline{Z}}s=0 \text{ for all } Z\in\Gamma(M,\mathcal{P})\},
\]
where $(\mathbb{L},h,\nabla)$ is a quantum bundle over $(M,E,\Omega)$.
Let us assume then that $M$ is a strongly pseudoconvex CR manifold with CR structure $E_{1,0}\subset T_\C M$.  Let $E\subset TM$ be the corresponding Levi distribution, and $J\in \End(E)$ the fibrewise complex structure on $E$ whose $+i$-eigenbundle is $E_{1,0}$. We choose a contact form (pseudo-Hermitian structure) $\theta$ such that the Webster metric
\[
 g_\theta(X,Y) = d\theta(JX,Y) + \theta(X)\theta(Y), \quad X,Y\in TM,
\]
is Riemannian.  Given this data, it is well-known (see \cite{DT}, for example) that there exists a unique linear connection $\nabla^{TW}$ on $M$, the {\em Tanaka-Webster connection}, such that:
\begin{enumerate}[(i)]
 \item $\nabla^{TW}_X \Gamma(M,E)\subset \Gamma(M,E)$, for all $X\in\Gamma(M,TM)$,
 \item $\nabla^{TW}J = \nabla^{TW}g_\theta = \nabla^{TW}\theta = 0$.
 \item The torsion $T_{TW}(X,Y)$ of $\nabla^{TW}$ is {\em pure}: for any $Z,W\in E_{1,0}$ and $X\in TM$, it satisfies
\begin{align*}
 T_{TW}(Z,W) & = 0\\
 T_{TW}(Z,\overline{W}) & = 2d\theta(Z,\overline{W})\xi\\
 T_{TW}(\xi, JX) & = - JT_{TW}(\xi,X),
\end{align*}
where $\xi$ denotes the Reeb field associated to $\theta$.
\end{enumerate}
Now, the Reeb field induces a splitting $T_\C M = E_{1,0}\oplus E_{0,1}\oplus \C\xi$ of the complexified tangent bundle.  Let us denote by $\widehat{T} = T_\C M/E_{0,1}\cong E_{1,0}\oplus\C\xi$.  We then obtain a bigrading of the space of complexified differential forms on $M$, given by
\[
 \mathcal{A}_\C^k(M) = \sum_{p+q = k}\mathcal{A}^{p,q}(M),
\]
where
\[
 \mathcal{A}^{p,q}(M) = \Gamma(M,\Lambda^p \widehat{T}^*\otimes\Lambda^q E^*_{0,1}).
\]
For $k = p+q$ we let $\pi^{p,q}:\mathcal{A}^k(M)\rightarrow \mathcal{A}^{p,q}(M)$ denote the corresponding projection.  Following \cite{BG} we define the {\em tangential Cauchy-Riemann operator}
\begin{equation}\label{gendbbar} 
\overline{\partial}_b = \pi^{p,q+1}\circ d:\mathcal{A}^{p,q}(M)\rightarrow \mathcal{A}^{p,q+1}(M),
\end{equation}
where $d:\mathcal{A}^r(M)\rightarrow \mathcal{A}^{r+1}(M)$ is the usual de Rham differential.
We note that for any $f\in C^\infty(M,\C)$ and $Z\in E_{1,0}$, we have
\[
 (\overline{\partial_b}f)(\overline{Z}) = \overline{Z}\cdot f.
\]
\begin{definition}
 We say that a function $f\in C^\infty(M,\C)$ is {\bf CR-holomorphic} if $\overline{\partial}_bf = 0$.
\end{definition}
We now introduce the CR analogue of a holomorphic vector bundle \cite{DT,Tanaka,Ura}:
\begin{definition}
Let $(M,E_{1,0})$ be a strongly pseudoconvex CR manifold. We say that a complex vector bundle $\mathcal{V}\rightarrow (M,E_{1,0})$ is {\bf CR-holomorphic} if it is equipped with a differential operator
\[
 \overline{\partial}_\mathcal{V}:\Gamma(M,\mathcal{V})\rightarrow \Gamma(M,E^*_{0,1}\otimes\mathcal{V})
\]
such that for any $u\in\Gamma(M,\mathcal{V})$, $f\in C^\infty(M,\C)$ and $Z,W\in E_{1,0}$,
\begin{align*}
 \overline{\partial}_\mathcal{V}(fu) & = f(\overline{\partial}_\mathcal{V} u) + (\overline{\partial}_b f)\otimes u\\
 [\overline{Z},\overline{W}]u &= \overline{Z}\,\overline{W}u = \overline{W}\,\overline{Z}u,
\end{align*}
where $\overline{Z}u = \iota(\overline{Z})(\overline{\partial}_\mathcal{V}u)$.
\end{definition}
Now, suppose we are given a CR-holomorphic vector bundle $\mathcal{V}\rightarrow (M,E_{1,0})$, equipped with a Hermitian metric $h$.  We say that a connection $\nabla$ on $\mathcal{V}$ is {\em Hermitian} if $\nabla h = 0$, and $\nabla^{0,1} := \nabla|_{E_{0,1}} = \overline{\partial}_\mathcal{V}$.  Such connections are uniquely determined up to a trace defined with respect to $\Omega = -d\theta$ \cite{DT,Ura}; the case where this trace is zero was introduced by Tanaka \cite{Tanaka}, and is known as Tanaka's canonical connection.

An example from \cite{DT} is the trivial line bundle $\mathbb{L} = M\times\C$ over $(M,E_{1,0})$, with the Hermitian metric $h_x((x,z_1),(x,z_2)) = z_1\overline{z}_2$.  The operator $\overline{\partial}_\mathbb{L}$ defined by
\begin{equation}\label{dbarl}
 (\overline{\partial}_\mathbb{L}s) = (x,(\overline{\partial}_b f)_x),\quad\text{for } s(x) = (x,f(x)),
\end{equation}
makes $\mathbb{L}$ into a CR-holomorphic line bundle.  If we equip $\mathbb{L}$ with the connection $\nabla$ defined by
\[
 \nabla_X s = (X\cdot f - i\theta(X))s,
\]
where $s(x) = (x,f(x))$, then the curvature form of $\nabla$ is equal to $i\Omega$, making $(\mathbb{L},h,\nabla)$ into a quantum bundle over $(M,E_{1,0},\Omega)$.  Moreover, the connection $\nabla$ is Hermitian; we have $\nabla^{0,1} = \overline{\partial}_\mathbb{L}$, and the following is therefore immediate:
\begin{proposition}\label{CRsectquant}
 Let $(M,E_{1,0})$ be a strongly pseudoconvex CR manifold with Levi distribution $E$. Let $\mathcal{P}$ be the CR polarization of $(M,E)$ given by $E_{1,0}$, and let $(\mathbb{L},h,\nabla)$ be the quantum bundle defined above.  Then the polarized sections of $\mathbb{L}$ are the {\em CR-holomorphic sections} of $\mathbb{L}$, defined by $\overline{\partial}_\mathbb{L}s = 0$.  Thus, $Q(M)$ is isomorphic to the space of CR-holomorphic $L^2$ functions on $M$.
\end{proposition}
\begin{remark}
This agrees with the answer given in \cite{BMG,GS2} for the homogeneous quantization of the symplectic cone given by the symplectization of an embedded strongly pseudoconvex CR manifold.
\end{remark}
\begin{remark}
In symplectic geometry, a polarization is given in general by a Lagrangian subbundle $\mathcal{P}$ of $TM\otimes\C$, with no condition on the rank of $\mathcal{P}\cap\overline{\mathcal{P}}$. (When this rank is maximal, the polarization is called {\em totally real}.)  It could be interesting to consider other versions of polarization in the contact setting.  The most natural definition of a totally real contact polarization should be given by the tangent bundle of a Legendrian foliation, but we have not yet considered this situation.
\end{remark}
\section{Dirac operators and index theory}
We will now review the index theoretic approach to the quantization of symplectic manifolds, before presenting an analogous theory for contact manifolds.  Numerous references now exist on this topic; a good overview can be found in \cite{SJ}.  We will refer primarily to the texts \cite{GGK,BGV}.
\subsection{The symplectic case}
We have already noted that given a prequanitzable K\"ahler manifold, a good candidate for the quantization $Q(M)$ is given by the space of holomorphic sections of the prequantum line bundle $\mathbb{L}$ (or some suitably high tensor power; see \cite{Mein1} for example).  When $(M,\omega)$ is a symplectic manifold, but not necessarily K\"ahler, we can mimic the above quantization with the aid of a suitable Spin$^c$-Dirac operator.

Given an even-dimensional Riemannian manifold $M$ with metric $g$, we can form the {\em Clifford bundle} $\Cl(M)\rightarrow M$ whose fibre over $x\in M$ is the complexified Clifford algebra of $T^*_xM$ with respect to the metric $g$. We let $\nabla$ denote a metric connection on $M$ (that is, $\nabla g=0$, but $\nabla$ is not necessarily torsion-free).  Such a connection preserves the Clifford multiplication, and induces a connection on the Clifford bundle $\Cl(M)$.
\begin{definition}
 A $\Z_2$-graded vector bundle $\mathcal{V}\rightarrow M$ is called a {\bf Clifford module} if there exists a homomorphism of graded algebras $a\in\Cl(M)\mapsto \mathbf{c}(a)\in\End(\mathcal{V})$.  We call a $\Z_2$-graded vector bundle $\mathcal{S}\rightarrow M$ a {\bf spinor bundle} if $\mathcal{S}$ is a Clifford module, and $\Cl(M)\rightarrow \End(\mathcal{S})$ is an isomorphism.
\end{definition}
\begin{remark}
Given any vector bundle $\mathcal{W}\rightarrow M$ and a Clifford module $\mathcal{V}\rightarrow M$, the tensor product bundle $\mathcal{V}\otimes\mathcal{W}$ is again a Clifford module, with respect to the Clifford action $\mathbf{c}(a)\otimes 1$.  If $M$ is equipped with a spin structure, then there is a canonical spinor bundle $\mathcal{S}$, and any other Clifford module is of the form $\mathcal{S}\otimes \mathcal{W}$ for some vector bundle $\mathcal{W}$.
\end{remark}
\begin{definition}
Let $\nabla$ be a metric connection on $(M,g)$. We say that a connection $\nabla^\mathcal{V}$ on a Clifford module $\mathcal{V}$ is a {\bf Clifford connection} (with respect to $\nabla$) if for every $a\in\Gamma(M,\Cl(M))$ and $X\in\Gamma(M,TM)$, we have
\[
 [\nabla^{\mathcal{V}}_X,\mathbf{c}(a)] = \mathbf{c}(\nabla_X a).
\]
\end{definition}
Given such a connection, we can define the geometric Dirac operator $\dob:\Gamma(M,\mathcal{V}^+)\rightarrow \Gamma(M,\mathcal{V}^-)$ given by $\dob = \mathbf{c}\circ\nabla^\mathcal{V}$ (see \cite{Nic} for example).
\begin{remark}
 Let $\mathcal{W}$ be a complex vector bundle with connection $\nabla^\mathcal{W}$, and let $\mathcal{S}$ be a spinor bundle equipped with a Clifford connection $\nabla^\mathcal{S}$.  The tensor product connection, given for $s\in\Gamma(M,\mathcal{S})$ and $w\in\Gamma(M,\mathcal{W})$ by
\begin{equation}\label{prodcon}
 \nabla^{\S\otimes\W}(s\otimes w) = \nabla^\mathcal{S}s\otimes w + s\otimes \nabla^\mathcal{W}w
\end{equation}
is then a Clifford connection on $\mathcal{S}\otimes \mathcal{W}$ with respect to the Clifford action $\mathbf{c}(a)\otimes\Id_\W$.
%One checks that
%\begin{equation*}
% [\nabla_X,\mathbf{c}(a)\otimes 1](s\otimes w) = [(\mathbf{c}(\nabla^{LC}_X a)\otimes 1](s\otimes w).
%\end{equation*}
\end{remark}
Let us now assume that $(M,\omega)$ is a symplectic manifold. Let $J\in \End(TM)$ be a compatible almost complex structure, and let $g$ be the corresponding Riemannian metric.  We can then choose a Hermitian connection $\nabla$ (such that $\nabla g = \nabla J = 0$).  Since $\nabla$ preserves the almost complex structure, it induces a connection $\nabla^\S$ on the spinor bundle
\[
 \mathcal{S} = \Lambda T^{0,1}M^*.
\]
The Clifford bundle $\Cl(E)$ acts on $\S$ via the action given by 
\begin{equation}\label{clifact}
 \mathbf{c}(\alpha)\nu = \sqrt{2}(\varepsilon(\alpha^{0,1}) - \iota(\alpha^{1,0}))\nu,
\end{equation}
where we have identified $T^{1,0}M^*$ with $T^{0,1}M$ by means of the Hermitian metric induced by $g$.  The induced connection $\nabla^\S$ is then a Clifford connection on $\mathcal{S}$.  If we assume that $(M,\omega)$ is prequantizable, and let $(\mathbb{L},\nabla^\mathbb{L},h)\to (M,\omega)$, then the tensor product connection on $\S\otimes\mathbb{L}$ is again a Clifford connection, and we can form the Dirac operator $\dob_\mathbb{L}$ acting on sections of $\S\otimes\mathbb{L}$.  The almost complex quantization of $(M,\omega)$ is then taken to be the $\Z_2$-graded vector space
\[
 Q(M) = \ker\dob_\mathbb{L}^+\oplus\ker\dob_\mathbb{L}^-,
\]
where $\dob_\mathbb{L}^+ = \dob_\mathbb{L}|_{\mathcal{A}^{0,2\bullet}(M,\mathcal{W})}$, and $\dob_\mathbb{L}^- =  \dob_\mathbb{L}|_{\mathcal{A}^{0,2\bullet + 1}(M,\mathcal{W})}=(\dob_\mathbb{L}^+)^*$.

Let us now specialize to the case where $(M,\omega)$ is K\"ahler.  In this case, the almost complex structure is integrable, and we can take the metric connection $\nabla$ to be the Levi-Civita connection $\nabla^{LC}$.
Given a Hermitian vector bundle $\mathcal{W}\rightarrow M$ with metric $h_\mathcal{W}$ and connection $\nabla^\mathcal{W}$, we have the decomposition $\nabla^\mathcal{W} = \nabla^{1,0} \oplus \nabla^{0,1}$ 
given by the restrictions of $\nabla^\mathcal{W}$ to $T^{1,0}M$ and $T^{0,1}M$, respectively.  When $\mathcal{W}$ is holomorphic, it is equipped with the differential operator
\[
 \overline{\partial}_\mathcal{W}:\mathcal{A}^{p,q}(M,\mathcal{W})\rightarrow \mathcal{A}^{p,q+1}(M,\mathcal{W})
\]
such that in any local holomorphic chart, $\overline{\partial}_\mathcal{W} = \sum\varepsilon(d\overline{z}^i)\dfrac{\partial}{\partial \overline{z}^i}$.  For $\mathcal{W}$ holomorphic, with Hermitian metric $h_\mathcal{W}$, then we have \cite[Proposition 3.65]{BGV}:
\begin{theorem}
 There exists a unique connection $\nabla^\mathcal{W}$ (the canonical connection) on $\mathcal{W}$ such that
\begin{enumerate}[(i)]
 \item $\nabla^{\mathcal{W}}h_\mathcal{W} = 0.$
 \item $\nabla^{0,1} = \overline{\partial}_\mathcal{W}$.
\end{enumerate}
\end{theorem}
Given a holomorphic vector bundle $\mathcal{W}$, the tensor product bundle $\mathcal{S}\otimes\mathcal{W}$ is a Clifford module, with the Clifford connection $\nabla$ given by \eqref{prodcon},
and we have the following \cite[Proposition 3.67]{BGV}:
\begin{theorem}
 The Dirac operator associated to the Clifford connection $\nabla^{\S\otimes\W}$ on $\mathcal{S}\otimes\mathcal{W}$ is given by
\begin{equation}\label{dob}
 \dob_\mathcal{W} = \sqrt{2}\left(\overline{\partial}_\mathcal{W} + \overline{\partial}_\mathcal{W}^{\,*}\right).
\end{equation}
\end{theorem}
The operator $\overline{\partial}_\mathcal{W}$ satisfies $\overline{\partial}_\mathcal{W}\circ \overline{\partial}_\mathcal{W} = 0$ on the complex
\[
 0\rightarrow \mathcal{A}^{0,0}(M,\mathcal{W})\rightarrow \mathcal{A}^{0,1}(M,\mathcal{W})\rightarrow \mathcal{A}^{0,2}(M,\mathcal{W})\rightarrow \cdots
\]
allowing us to define the Dolbeault cohomology groups $H^{0,q}(M,\mathcal{W})$, which in turn are isomorphic to the sheaf cohomology groups $H^q(M,\mathcal{O}(\mathcal{W}))$, where $\mathcal{O}(\mathcal{W})$ denotes the sheaf of holomorphic sections of $\mathcal{W}$.

By the usual Hodge theory for the Dolbeault complex, we obtain the equality of $\Z_2$-graded vector spaces
\[
 \ker(\dob^+_\mathcal{W})\oplus\ker(\dob^-_\mathcal{W}) = \sum(-1)^iH^{0,i}(M,\mathcal{W}) \cong \sum(-1)^iH^i(M,\mathcal{O}(\mathcal{W})).
\]
For the case $\mathcal{W} = \mathbb{L}$, if we assume that $\mathbb{L}$ is sufficiently positive, then $H^k(M,\mathcal{O}(\mathcal{W}))=0$ for $k>0$ by the Kodaira vanishing theorem, and the almost complex quantization of $(M,\omega)$ coincides with the earlier definition in terms of the space of holomorphic sections of $\mathbb{L}$.  Moreover, the dimension of $Q(M)$ is given by the index of the Dirac operator $\dob_\mathbb{L} = \sqrt{2}\left(\overline{\partial}_\mathbb{L}+\overline{\partial}_\mathbb{L}^*\right)$,  which we can compute via the Riemann-Roch formula:
\[
 \ind(\dob_\mathbb{L}) = \frac{1}{(2\pi i)^n}\int_M \Td(TM)\Ch(\mathbb{L}).
\]
When $M$ is equipped with a Hamiltonian $G$-action $Q(M)$ becomes a virtual $G$-representation, and the associated virtual character is given near the identity, for $X\in\mathfrak{g}$ sufficiently small, by the equivariant Riemann-Roch formula
\begin{equation}\label{RR}
\chi(e^X) = \ind^G(\dob_\mathbb{L})(e^X) = \frac{1}{(2\pi i)^n}\int_M \Td(TM,X)\Ch(\mathbb{L},X),
\end{equation}
with similar formulas near other elements $g\in G$ (see \cite{BGV}).
\subsection{The Contact Case}
Let us now assume that $(M,E)$ is a contact manifold, and let $\theta$ be a choice of contact form on $M$.  Since $\Omega = -d\theta|_{E\otimes E}$ makes $E\to M$ into a symplectic vector bundle, we can choose a compatible complex structure $J\in\End(E)$ on the fibres of $E$. The $+i$-eigenbundle $E_{1,0}\subset T_\C M$ of $J$ then defines an almost CR structure on $M$.  If we extend $J$ to $TM$ by setting $J\xi=0$, where $\xi$ is the Reeb vector field associated to $\theta$, then $J$ is an almost contact structure on $M$.  We let $g$ be a Riemannian metric on $M$ such that $(J,g,\theta,\xi)$ is a contact metric structure.  (Precisely, we can set $g(X,Y) = d\theta(JX,Y)+\theta(X)\theta(Y)$.) We let $\Cl(E)$ denote the bundle of Clifford algebras over $M$ whose fibre over $x\in M$ is the complexified Clifford algebra of $E^*_x$ with respect to Euclidean form given by restricting $g$ to $E$.

It is always possible to choose a connection $\nabla$ such that $\nabla\theta =\nabla J = \nabla g = 0$; such connections are known as contact metric connections in \cite{Nic}.  Since $\nabla\theta = 0$, it follows that $\nabla$ preserves the contact distribution $E$, and thus, since $\nabla g=0$, we obtain an induced connection on $\Cl(E)$ that is compatible with the Clifford multiplication in $\Cl(E)$.

As in the even-dimensional case, we will call a $\Z_2$-graded vector bundle $\mathcal{V}\rightarrow M$ a {\em Clifford module} if there exists a homomorphism of graded algebras $\Cl(E)\rightarrow \End(\mathcal{V})$, and we will call a $\Z_2$-graded vector bundle $\mathcal{S}\rightarrow M$ a {\em spinor bundle} if $\Cl(E)\rightarrow \End(\mathcal{V})$ is an isomorphism of graded algebras.

Using the almost CR structure $E_{1,0}$, we can define the bundle $\S = \Lambda E_{0,1}^*$.  Since $\nabla$ preserves $E$ and $\nabla J=0$, it induces a connection $\nabla^\S$ on $\mathcal{S}$.  The bundle $\S$ is then a spinor bundle for $\Cl(E)$, where we define the Clifford action of $\Cl(E)$ on $\mathcal{S}$ by \eqref{clifact}, keeping in mind that we must take $\alpha\in E^*$ and not $\alpha\in T^*M$. (Note that the splitting of $T^*M$ determined by the contact form $\theta$ allows us to identify $E^*$ with a subbundle of $T^*M$.)  By analogy with the symplectic case, we make the following definition:
\begin{definition}
Let $\nabla$ be a contact metric connection on $M$.  We say that a connection $\nabla^\mathcal{V}$ on a Clifford module $\mathcal{V}$ is a {\bf Clifford connection} with respect to $\nabla$ if for every $a\in\Gamma(M,\Cl(E))$ and $X\in\Gamma(M,TM)$, we have
\begin{equation}\label{CRclif}
 [\nabla^{\mathcal{V}}_X,\mathbf{c}(a)] = \mathbf{c}(\nabla_X a).
\end{equation}
\end{definition}
\begin{proposition}
The connection $\nabla^\S$ induced on $\mathcal{S} = \Lambda E^*_{0,1}$ by the contact metric connection $\nabla$ is a Clifford connection.
\end{proposition}
\begin{proof}
Since the $\nabla$ is compatible with the Clifford multiplication, it suffices to check that \eqref{CRclif} holds for a 1-form $\alpha\in\Gamma(M,E^*)$.  For any $s\in\Gamma(M,\mathcal{S})$, we have
\begin{align*}
[\nabla^\S_X,\mathbf{c}(\alpha)]s &= (\nabla_X\alpha^{0,1})\wedge s +\alpha^{0,1}\wedge\nabla^\S_X s\\
			       &\hspace{40pt}-(\iota(\nabla_X\alpha^{1,0})s-\iota(\alpha^{1,0})\nabla^\S_X s)\\
			       &\hspace{80pt}-(\alpha^{0,1}\wedge\nabla^\S_X s -\iota(\alpha^{1,0})\nabla^\S_X s)\\
			       &= (\nabla_X\alpha^{0,1})\wedge s -\iota(\nabla_X\alpha^{1,0})s\\
			       &= \mathbf{c}(\nabla_X\alpha)s.
\end{align*}
\end{proof}
Given such a Clifford connection on a Clifford module $\mathcal{V}=\mathcal{V}^+\oplus\mathcal{V}^-$, we can define a Dirac-like operator $\dirac:\Gamma(M,\mathcal{V}^+)\rightarrow\Gamma(M,\mathcal{V}^-)$ by the composition
\begin{equation}\label{dirac}
 \dirac: \Gamma(M,\mathcal{V}^+)\xrightarrow[]{\nabla^\mathcal{V}}\Gamma(M,T^*M\otimes\mathcal{V}^+)\xrightarrow[]{q}\Gamma(M,E^*\otimes\mathcal{V}^+)\xrightarrow[]{\mathbf{c}}\Gamma(M,\mathcal{V}^-),
\end{equation}
where $q:T^*M\rightarrow E^*$ is orthogonal projection with respect to the metric $g$.  Given an auxilliary complex vector bundle $\mathcal{W}\rightarrow M$ with connection $\nabla^\mathcal{W}$, the tensor product connection $\nabla^{\S\otimes\W}$ given by \eqref{prodcon} is again a Clifford connection with respect to $\nabla$ on $\mathcal{S}\otimes\mathcal{W}$: for any section $s\otimes w\in\Gamma(M,\mathcal{S}\otimes \mathcal{W})$, we have
\begin{align*}
 [\nabla^{\S\otimes\W}_X,\mathbf{c}(a)\otimes 1]s\otimes w & = \nabla^{\S\otimes\W}_X (\mathbf{c}(a) s\otimes w) - \mathbf{c}(a)\otimes 1(\nabla^\mathcal{S}_X s\otimes w+s\otimes\nabla^\mathcal{W}_X w)\\
&=([\nabla_X,\mathbf{c}(a)] s)\otimes w\\
&=(\mathbf{c}(\nabla_X a)s)\otimes w\\
&=(\mathbf{c}(\nabla_X a)\otimes 1)s\otimes w.
\end{align*}
Thus, we can define the twisted Dirac operator $\dob_\W$ acting on sections of $\S\otimes\W$.  In particular, let us consider the trivial bundle $\mathbb{L}=M\times\C\rightarrow M$, equipped with the connection $\nabla^\mathbb{L}$ and Hermitian metric $h$ defined in Section \ref{crpolar}.  Then $(\mathbb{L},\nabla^\mathbb{L},h)\to (M,E,\Omega)$ is again a quantum bundle (but not yet a CR-holomorphic line bundle, since we are not assuming that $M$ is CR at the moment).  Using the tensor product connection on $\mathcal{S}\otimes\mathbb{L}$, we can define the Dirac operator 
\[
\dob^\pm_\mathbb{L}:\Gamma(M,\mathcal{S}^\pm\otimes\mathbb{L})\to\Gamma(M,\mathcal{S}^\mp\otimes\mathbb{L}),
\]
and define the ``almost CR quantization'' of our contact manifold $M$ as
\begin{equation}\label{acrquant}
Q(M) = \ker(\dob^+_\mathbb{L})\oplus\ker(\dob^-_\mathbb{L}).
\end{equation}
Now, unlike in the symplectic case, the operator $\dob_\mathbb{L}$ is not elliptic.  This is not entirely bad, since on an odd-dimensional manifold the index of an elliptic operator is trivial.  On the other hand, given a $G$-action commuting with $\dob_\mathbb{L}$, the virtual representation given by $Q(M)$ will be infinite-dimensional, so in general we do not have an analogue of the character forumula given by \eqref{RR}.

However, as shown in \cite{F}, when $G$ acts on $M$ such that the orbits of the $G$-action are transverse to the contact distribution $E$, $\dob_\mathbb{L}$ is a transversally elliptic operator.  Thus, by a result of Atiyah and Singer \cite{AT}, the equivariant index of $\dob_\mathbb{L}$ is well-defined as a generalized function on $G$.  Using the contact form $\theta$ we can define an equivariant differential form with generalized coefficients $\mathcal{J}(E,X)$ (that is, a differential form depending distributionally on $X\in\mathfrak{g}$ \cite{KV}) given by
\begin{equation}\label{jex}
 \mathcal{J}(E,X) = \theta\wedge\delta_0(D\theta(X)),
\end{equation}
where $\delta_0(x)$ denotes the Dirac delta distribution on $\R$, and $D\theta(X)=d\theta - \Phi^X_\theta$ is the equivariant differential of $\theta$.  This form was introduced in \cite{F}; it is equivariantly closed, and depends only on the contact distribution $E$ and the group action.  In particular, it is independent of the choice of contact form $\theta$. By manipulating the equivariant index formulas of Berline-Paradan-Vergne \cite{BV1,BV2,PV3}, we showed that the equivariant index of the $G$-transversally elliptic operator $\dirac$ is then given near $1\in G$ by
\[
 \ind^G(\dob_\mathbb{L})(e^X) = \frac{1}{(2\pi i)^n}\int_M \Td(E,X)\Ch(\mathbb{L},X)\mathcal{J}(E,X),
\]
for $X\in\mathfrak{g}$ sufficiently small. The formula near other elements $g\in G$ is similar; the integration is then over the fixed-point set $M^g$, and the integrand includes a contribution from the normal bundle.  (We showed in \cite{F} that the corresponding restriction of $\mathcal{J}(E,X)$ is well-defined.)   One advantage of defining $Q(M)$ in terms of the equivariant index of $\dob_\mathbb{L}$ is that it follows that $Q(M)$ does not depend on the choice of contact form $\theta$, whereas the geometric construction involves a choice of contact form throughout.

Now, we would like to be able to relate the above index-theoretic case to the quantization defined earlier for the case of a strongly pseudconvex CR manifold.  Let us assume then that the subbundle $E_{1,0}$ is involutive, so that it determines a CR structure on $M$ for which the contact form $\theta$ is a pseudo-Hermitian structure.  In this case the metric $g$ is the Webster metric associated to $\theta$, and there is a canonical contact metric connection; namely, the Tanaka-Webster connection $\nabla^{TW}$.

By \cite[Proposition 1.17]{DT}, the $\overline{\partial}_b$ operator can be written in terms of a local frame $\{Z_i\}$ for $E_{1,0}$ with corresponding coframe $\{\theta^i\}$ according to
\[
 \overline{\partial}_b \alpha= \sum \overline{\theta}^i\wedge (\nabla^{TW}_{\overline{Z}_i}\alpha),
\]
for any $\alpha\in\mathcal{A}^{0,k}(M)$.  In other words, as an operator on $\mathcal{A}^{0,\bullet}(M)$, $\overline{\partial}_b$ is given by the composition
\[
 \mathcal{A}^{0,k}(M)\xrightarrow{\nabla^{TW}}\Gamma(M,T^*M\otimes \Lambda^kE^*_{0,1})\xrightarrow{q}\Gamma(M,E^*\otimes\Lambda^k E^*_{0,1})\xrightarrow{\varepsilon}\mathcal{A}^{0,k+1}(M),
\]
where $q:T^*M\rightarrow E^*$ is orthogonal projection with respect to $g$, and for any $\alpha\in \Gamma(M,E^*)$ and $\gamma\in\mathcal{A}^{0,k}(M)$, we define $\varepsilon(\alpha)\cdot \gamma =\alpha^{0,1}\wedge\gamma$.

Some care must be taken in obtaining the above decomposition of the $\overline{\partial}_b$ operator: since $\nabla^{TW}$ has torsion, we cannot write the full exterior differential $d$ in terms of the Tanaka-Webster connection.  The proof relies on the fact that the torsion of $\nabla^{TW}$ is pure, and hence vanishes on $E_{0,1}\otimes E_{0,1}$.

Similarly, (see equation (1.142) in \cite{DT}, and the line immediately above it) the formal adjoint of $\overline{\partial}_b$ is given locally by the expression $\overline{\partial}_b^*\gamma = -\sum \iota(\overline{Z}_i)(\nabla_{Z_i}\gamma)$.  Globally, we write this as the composition
\[
\mathcal{A}^{0,k}(M)\xrightarrow{\nabla^{TW}}\Gamma(M,T^*M\otimes\Lambda^k E^{0,1})\xrightarrow{q}\Gamma(M,E^*\otimes\Lambda^k E^{0,1})\xrightarrow{-\iota}\mathcal{A}^{k-1}(M),
\]
where $\iota(\alpha)\cdot \gamma = \iota(\alpha^{1,0})\gamma$.  Here $\alpha^{1,0}\in E^*_{1,0}$, and we identify $E^*_{1,0}=\overline{E^*_{0,1}}=E_{0,1}$ using the Hermitian metric determined by $g_\theta$.  Thus, we obtain the following:
\begin{proposition}
 The Dirac operator $\dirac$ associated to the connection $\nabla^\S$ induced by the Tanaka-Webster connection $\nabla^{TW}$  on $\mathcal{S}=\Lambda E^*_{0,1}$ is given by
\[
 \dirac = \sqrt{2}\left(\overline{\partial}_b + \overline{\partial}_b^*\right).
\]
\end{proposition}
We can thus interpret the index of $\dirac$ in terms of the Kohn-Rossi cohomology groups \cite{KR}
\begin{equation}
 H^{0,i}(M,E_{1,0}) = \frac{\ker(\overline{\partial}_b:\mathcal{A}^{0,i}(M)\rightarrow \mathcal{A}^{0,i+1}(M))}{\im(\overline{\partial}_b:\mathcal{A}^{0,i-1}(M)\rightarrow \mathcal{A}^{0,i}(M))}.
\end{equation}
From \cite{Kohn}, we have that $H^{0,i}(M,E_{1,0}) \cong \ker(\Box^i_b)$, where $\Box^i_b$ is the Kohn-Rossi Laplacian acting on $(0,i)$-forms. Since $\ker(\Box_b) = \ker(\overline{\partial}_b)\cap \ker(\overline{\partial}_b^*)$ \cite{Kohn} and $\Box_b=\dirac^{\!2}$, we have $\ker(\Box_b^i) = \ker(\dirac^{\! i})$. Thus, we may prefer to use the equivalent definition of $Q(M)$ given by
\begin{equation}\label{newqm}
Q(M) = \sum_{i=0}^n (-1)^iH^{0,i}(M,E_{1,0}).
\end{equation}
However, this is not entirely satisfactory, since $\ker\dirac$ is in general infinite-dimensional. (Kohn proved in \cite{Kohn} that $H^{0,i}(M,E_{1,0})$ is finite-dimensional, but only for $1\leq i\leq n-1$.)  When $Q(M)$ is infinite-dimensional, we should expect it to be a Hilbert space.  Thus, it may make more sense to replace the Kohn-Rossi cohomology groups in the above definition of $Q(M)$ by suitable $L^2$ cohomology groups for the $\dbbar$ operator, and since we are restricting ourselves to strongly pseudoconvex CR manifolds, we have a preferred metric $g_\theta$ and measure $\mu_\theta$ with which to define the $L^2$ cohomology.  We have not yet investigated the consequences of using $L^2$ cohomology groups to define $Q(M)$, but we note that one advantage of doing so is that in degree zero, we recover our earlier definition of $Q(M)$ in terms of the CR-holomorphic $L^2$ functions on $M$.  However, we are not aware of any analogue of the Kodaira vanishing theorem that would allow us to identify our two versions of contact quantization in certain settings.
\begin{remark}
For a brief discussion of the use of $L^2$ cohomology groups for non-compact symplectic manifolds, see \cite{GGK}.  In particular, we note \cite[Remark 6.36]{GGK}, which points out that in general, defining a virtual representation of the form \eqref{newqm} may not make sense in general, when some of the terms are infinite-dimensional.  However, they do note that such an expression is well-defined as a virtual representation whenever every irreducible representation occurs with finite multiplicity, which is in particular the case whenever $\dirac$ is a transversally elliptic operator.
\end{remark}

Finally, we note that if we are given a CR-holomorphic vector bundle $(\mathcal{W},\overline{\partial}_\mathcal{W})\rightarrow M$, equipped with a Hermitian metric $h_\mathcal{W}$, we can assume that $\mathcal{W}$ is equipped with a Hermitian connection $\nabla^\mathcal{W}$.  The bundle $\mathcal{S}\otimes\mathcal{W}$ will then be a Clifford bundle, with Clifford connection $\nabla^{\S\otimes\W}$ given by the tensor product connection \eqref{prodcon}. Moreover, the operator $\overline{\partial}_\W$ extends to an operator
\[
\overline{\partial}_\W:\mathcal{A}^{0,k}(M,\W)\to\mathcal{A}^{0,k+1}(M,\W)
\]
on $\W$-valued $(0,k)$-forms such that $\overline{\partial}_\W^2=0$ \cite{DT}. We then have the following:
\begin{theorem}
 The Dirac operator on $\mathcal{S}\otimes\mathcal{W}$ determined by the Clifford connection $\nabla^{\S\otimes\W}$ is given by
\[
 \dob_\mathcal{W} = \sqrt{2}\left(\overline{\partial}_\mathcal{W} + \overline{\partial}_\mathcal{W}^{\,*}\right),
\]
where $\overline{\partial}_\mathcal{W}$ is the extension of the CR holomorphic operator of $\mathcal{W}$ to $\mathcal{S}\otimes\mathcal{W}$.
\end{theorem}
\begin{proof}
 We need to show that the two operators agree on sections of $\mathcal{W}\otimes\Lambda E^*_{0,1}$.  Let $\{\overline{Z}_i\}$ be a local frame for $E_{0,1}$, with corresponding coframe $\{\overline{\theta}^i\}$ for $E^*_{0,1}$.  We note that $\overline{\partial}_\mathcal{W}$ can be expressed locally by
\begin{equation}\label{dwloc}
 \overline{\partial}_\mathcal{W}\alpha = \sum \overline{\theta}^i\wedge \overline{Z}_i\alpha,
\end{equation}
where $\overline{Z}_i\alpha = \iota(\overline{Z}_i)\overline{\partial}_\mathcal{W}\alpha$.  Let $\nabla$ be the tensor product connection on $\Lambda E^*_{0,1}\otimes \mathcal{W}$, and define the operator $\overline{\partial}_\nabla:\Gamma(M,\Lambda^k E^*_{0,1}\otimes\mathcal{W})\rightarrow \Gamma(M,\Lambda^k E^*_{0,1}\otimes\mathcal{W})$ given by
\[
 \overline{\partial}_\nabla \alpha = \sum \overline{\theta}^i\wedge (\nabla_{\overline{Z}_i}\alpha).
\]
From \cite{Ura}, we have that $\overline{\partial}^*_\nabla = -\sum\iota(\overline{Z}_i)\nabla_{Z_i}$. It follows that $\dob_\mathcal{W} = \sqrt{2}\left(\overline{\partial}_\nabla + \overline{\partial}_\nabla^*\right)$, so it suffices to show that as operators on $\Lambda E^*_{0,1}\otimes\mathcal{W}$, we have $\overline{\partial}_\nabla = \overline{\partial}_\mathcal{W}$.  Let $\alpha\otimes w\in\Gamma(M,\Lambda E^*_{0,1}\otimes\mathcal{W})$.  Then, since $\nabla^\mathcal{W}$ is a Hermitian connection, we have that for any $\overline{Z}\in E_{0,1}$, $\nabla^\mathcal{W}_{\overline{Z}}w = \iota(\overline{Z})(\overline{\partial}_\mathcal{W}w)$, and therefore,
\begin{align*}
 \overline{\partial}_\nabla(\alpha\otimes w) &= \sum \overline{\theta}^i\wedge\left(\nabla^{TW}_{\overline{Z}_i}\alpha\otimes w + \alpha\otimes \nabla^\mathcal{W}_{\overline{Z}_i}w\right)\\
&= \sum\left(\overline{\theta}\wedge(\nabla^{TW}_{\overline{Z}_i}\alpha)\otimes w + \overline{\theta}^i\wedge\alpha\otimes (\iota(\overline{Z}_i)\overline{\partial}_\mathcal{W}w)\right)\\
&= (\overline{\partial}_b\alpha)\otimes w + (-1)^{|\alpha|}\alpha\wedge \sum \overline{\theta}^i\otimes \overline{Z}_iw\\
&= (\overline{\partial}_b\alpha)\otimes w + (-1)^{|\alpha|}\alpha\wedge (\overline{\partial}_\mathcal{W} w)\\
&= \overline{\partial}_\mathcal{W}(\alpha\otimes w).\qedhere
\end{align*}
\end{proof}
\bibliographystyle{amsalpha}
\bibliography{references}

\end{document}